\documentstyle[11pt, leqno]{article}

\begin{document}
\def\s{\subseteq}
\def\n{\noindent}
\def\se{\setminus}
\def\dia{\diamondsuit}
\def\la{\langle}
\def\ra{\rangle}


\title{Multiplicative Zagreb indices of cacti}

\footnotetext{The first author was
partially supported by the Summer Graduate Research Assistantship Program of Graduate School,
the second author was partially supported by College of Liberal Arts Summer Research Grant. }

\author{Shaohui Wang\footnote{  Corresponding author: S. Wang (e-mail: swang4@go.olemiss.edu), B. Wei (e-mail:  bwei@olemiss.edu). } , Bing Wei \\ 
\small\emph {Department of Mathematics, The University of Mississippi, University, MS 38677, USA}\\
Accepted by Discrete Mathematics, Algorithms and Applications}
\date{}
\maketitle

\begin{abstract}

Let   $\prod(G)$ be  Multiplicative Zagreb index of  a graph $G$.   A connected graph is a cactus graph if and only if any two of its cycles have at most one vertex in common, which has been the interest of researchers in the filed of material chemistry and graph theory.   In this paper, we use a new tool  to the obtain upper and lower bounds of $\prod(G)$  for all cactus graphs and  characterize the corresponding extremal graphs.

\end{abstract}

\section{Introduction}
 During recent decades, applied graph theory, molecular topology and mathematical chemistry have been the focus of considerable research  in developed theory.   In the field of  chemical molecular graphs [9,14], the atoms are represented by vertices and the bonds by edges that capture the structural essence of compounds.  
The numerical  representation of the molecule graph
 can be mathematically deduced as a single number, usually called graph invariant, molecular descriptor or  topological index.  

One of the oldest and most thoroughly considered  molecular descriptor is  Zagreb index which was introduced by
  Gutman and Trinajsti$\acute{c}$ in 1972[6] as below:
The first Zagreb index $M_1$ is the sum of the square vertex degrees of all the atoms and the second Zagreb index $M_2$ is the sum over all  bonds of the product of the vertex degrees of the two  adjacent atoms, that is, for any graph $G = (V, E)$  with vertex set $ V(G)$ and edge set $E(G)$,
 $$ M_1(G) = \sum_{u \in V(G)} d(u)^2,
M_2(G) = \sum_{uv \in E(G)} d(u)d(v).$$ 
In the 1980s, Narumi and Katayama [14] characterized  the structural isomers of saturated hydrocarbons and considered the product 
$
NK = \prod_{v \in V(G)} d(v),
$
which  is called  the "Narumi-Katayama index". Recently, Todeschini et al.[5,7],  Wang and Wei[16]  studied the first (generalized) and second Multiplicative Zagreb indices defined as follows: For $c > 0$, 
$$
\prod_1(G) = \prod_{v \in V(G)}d(v)^{c},
\prod_2(G) = \prod_{uv \in E(G)} d(u)d(v).
$$
Obviously, the first Multiplicative Zagreb index is the power of the NK index.  Moreover, the second Multiplicative Zagreb index can be rewritten as $\prod_2(G) = \prod_{v \in V(G)}d(v)^{d(v)}$.

In the past several years, there are a lot of significant and interesting results [8,15] about chemical indices to the study of a computational complexity and the intersection between graph theory and chemistry.  
For general graphs, a lower bound of a chemical index, called Randi$\acute{c}$ index, was given by Bollobás and Erdös(1998)[1], while an upper bound was
recently presented in (2004)[13].
In 2004, Das [2] applied the minimal and maximal degree to obtain the upper bound for the sum of the squares of the degrees of a
graph, the first Zagreb index.  In 2010, Zhao and Li [18] provided
the maximal Zagreb index of graphs with $k$ cut vertices.
 Estes and Wei (2014)[4], Wang and Wei (2015)[16] gave the sharp upper and lower bounds of  Zagreb indices and Multiplicative Zagreb indices of $k$-trees, a generalization of a tree, respectively. 

The synthetic resins\cite{pl}, a type of plastic materials, is produced by the condensation of phenol with
formaldehyde in the presence of a base.  Independent benzene rings have no common edges in the diphenyl ether and the biphenyl. Many generic phenolic  structures of clindamycin phosphate and  cellulose  have no shared edges  between different phenolic molecules.  Due to the properties, we consider a special class of graphs:
A graph is a cactus if it is connected and all of its blocks  are either
edges or cycles, i.e., any two of its cycles have at most one common vertex.   In 1969, Cornu$\acute{e}$jols and Pulleyblank [3] used the constructure of a triangular cactus  to find the equivalent conditions for the existence of $\{K_2, C_n, n \geq 4\}$-factor. Also, 
 Lin et al.(2007)[11],  Liu and Lu (2008)[12]  obtained some sharp bounds of several chemical indices of cactus graphs, such as Wiener index, Merrifield-Simmons index, Hosoya index and
Randi$\acute{c}$ index. In 2012, Li et al.[10] gave the upper  bounds on Zagreb  indices of  cactus graphs and  lower bounds of cactus graph with at least one cycle.   Wang and Kang(2015) [17] found  the extremal   bounds of another chemical index, Harary index, for the cacti as well.

In this paper, we  use the new tool of  interesting functions to  obtain  sharp bounds of Multiplicative Zagreb indices on cacti, which can partially indicate the strength of heat resistance and flame retardancy by maximal and minimal bounds.  By taking the derivatives, one can check the following facts.
\vskip 2mm {\bf  Fact 1. }  \emph {
The function $f(x) =\displaystyle{\frac {x}{x+m}}$ is strictly increasing for  $x \in [0, \infty)$, where $m$ is a positive integer.
}
\vskip 2mm {\bf  Fact 2. }  \emph {
The function $f(x) =\displaystyle{\frac {x^x}{(x+m)^{x+m}}}$  is strictly decreasing for  $x \in [0, \infty)$, where $m$ is a positive integer.
}

  Since every cactus graph may have some pendant vertices which connect to one vertex  only, then set $\mathcal{C}$$_n^k$ to denote a set of cactus graphs with $n$ vertices including $k$ pendant vertices, where $n \geq k \geq 0$. An edge is called a  pendant  edge if  one of its vertices is a pendant vertex.
For $r \geq 1$, let $P_1 = u_1u_2...u_pv_1, P_2 = u_1u_2...u_pv_2, ..., P_r = u_1u_2...u_pv_r$ be the paths of a graph $G$ such that there exists at most one cycle $C$ with $ V(P_i) \cap V(C) = \{u_1\}$ and $d(v_i) = 1$, $i \geq 1$,  then  the induced subgraph $G[\{v_i, u_j, i \in [1,r], j \in [1,p]\}]$ is called a dense path. In particular, when $r =1$, the dense path is a pendant path. The length of a dense path is the length of its  pendant path. 
Theorems 1,2,3  provide   sharp upper and lower  bounds on the first generalized Multiplicative Zagreb  indices of  cactus graphs  and characterize the  extremal graphs. 
\vskip 2mm {\bf Theorem 1 } \emph{
For any graph $G $ in $\mathcal{C}$$_n^k$, 
$$\prod_{1,c}(G) \geq \left\{ \begin{array}{rcl}
3^{kc}2^{(n-2k)c} & \mbox{if} &  k=0, 1, \\ 
2^{(n-k-1)c} k^c & \mbox{if} &  k \geq 2,
\end{array}\right.$$
 the equalities hold if and only if  their degree sequences are $\underbrace{3, 3, ..., 3}_k, \underbrace{2,2, ..., 2}_{n-2k}, \underbrace{1,1,...,1}_k$ and   $k,  \underbrace{2,2,  ..., 2}_{n-k-1}, \\ \underbrace{1,1, ..., 1}_{k}$, respectively.
}

\vskip 2mm {\bf Theorem 2 } \emph{
For any graph $G $ in $\mathcal{C}$$_n^k$ with $ n \leq k+3$,
$$\prod_{1,c}(G) \leq \left\{ \begin{array}{rcl}
 k^c & \mbox{if} &  n=k+1,\\ 
 (\lceil{\frac{k}{2}}\rceil+1)^c(\lfloor{\frac{k}{2}}\rfloor+1)^c & \mbox{if} & n= k+2,\\ 
(\lceil{\frac{k}{3}}\rceil+2)^c(\lfloor{\frac{k}{3}}\rfloor+2)^c(k-\lceil{\frac{k}{3}}\rceil-\lfloor{\frac{k}{3}}\rfloor+2)^c & \mbox{if} &  n=k+3,\\ 
\end{array}\right.$$
 the equalities hold if and only if  their degree sequences  are   $k,\underbrace{1,1,...,1}_{k}$; $\; \lceil{\frac{k}{2}}\rceil+1, \lfloor{\frac{k}{2}}\rfloor+1, \underbrace{1,1,...,1}_{k}$ and  $\lceil{\frac{k}{3}}\rceil+2,\lfloor{\frac{k}{3}}\rfloor+2, k-\lceil{\frac{k}{3}}\rceil-\lfloor{\frac{k}{3}}\rfloor+2,\underbrace{1,1,...,1}_{k},$  respectively.
}

\vskip 2mm {\bf Theorem 3 } \emph{
For any graph $G $ in $\mathcal{C}$$_n^k$ with $ n \geq k+4$  and $t \geq 0$, 
$$\prod_{1,c}(G) \leq \left\{ \begin{array}{rcl}
16^c  & \mbox{if} &  k = 0, n=4,\\ 
2^{(3t+6)c} & \mbox{if} &  k = 0, n=2t+5,\\ 
2^{(3t+4)c}9^c  & \mbox{if} & k = 0, n= 2(t+3),\\ 
\end{array}\right.$$
 the equalities hold if and only if  their degree sequences are $2,2,2,2$;  $ \underbrace{4,4,...,4}_{t+1}, \underbrace{2,2,...,2}_{t+4}$ and  $\underbrace{4,4,...,4}_{t},   3, 3,\\  \underbrace{2,2,...,2}_{t+4}$,  respectively; \\ For $k \neq 0$, if $\prod_{1,c}(G)$ attains the maximal value, then one of the following statements holds: For any nonpendant vertices $u,v$, either (i) $|d(u) - d(v)| \leq 1$ or (ii) $d(u) \in  \{2,3,4\}$ and $G$ contains no cycles of length greater than 3, no dense paths of length greater than 1 except for at most one of them with length 2,  no paths of length greater 0  that connects only two cycles except for at most one of them with length 1. 
}

Theorems 4, 5  give  the sharp upper and lower bounds on the second Multiplicative Zagreb  indices of  cactus graphs and characterize the extremal graphs.
\vskip 2mm {\bf Theorem 4 } \emph{
For any graph $G $ in $\mathcal{C}$$_n^k$ with $\gamma = \frac{k-2}{n-k}$, 
$$\prod_{2}(G) \geq \left\{ \begin{array}{rcl}
3^{3k}2^{2(n-2k)} & \mbox{if} &  k =0, 1, \\ 
(2 +\lceil{ \gamma}\rceil )^{(2 + \lceil{ \gamma}\rceil )[k-2 - \lfloor{ \gamma}\rfloor(n-k)]} 
(2 +  \lfloor{\gamma}\rfloor) ^{(2 +\lfloor{ \gamma}\rfloor )[n- 2k+2 + \lfloor{ \gamma}\rfloor(n-k)]}
& \mbox{if} & k\geq 2,
\end{array}\right.$$
the equalities hold if and only if  their degree sequences  are  $\underbrace{3,3, ...,3}_k, \underbrace{2,2,  ..., 2}_{n-2k}, \underbrace{1,1,...,1}_k$ and\\ $\underbrace{2+\lceil{\gamma}\rceil, 2+ \lceil{\gamma}\rceil, ..., 2+ \lceil{\gamma}\rceil}_{k-2 -\lfloor{\gamma}\rfloor (n-k) }, \underbrace{2+ \lfloor{\gamma}\rfloor, 2+ \lfloor{\gamma}\rfloor, ..., 2+\lfloor{\gamma}\rfloor}_{n-2k+2+ \lfloor{\gamma}\rfloor(n-k)}, \underbrace{1, 1, ..., 1}_{k}$, respectively.
}

\vskip 2mm {\bf Theorem 5 } \emph{
For any graph $G $ in $\mathcal{C}$$_n^k$, 
$$\prod_{2}(G) \leq \left\{ \begin{array}{rcl}
(n-2)^{n-2}2^{2(n-k-1)} & \mbox{if} &  n-k \equiv 0 (\mbox{mod}\; 2), \\ 
(n-1)^{n-1}2^{2(n-k-1)} & \mbox{if} &  n-k \equiv 1 (\mbox{mod}\; 2),
\end{array}\right.$$
the equalities hold if and only if  their degree sequences   are $n-2, \underbrace{2,2,  ..., 2}_{n-k-1 }, \underbrace{1, 1, ...,1}_{k}$ and $n-1, \underbrace{2, 2, ..., 2}_{n-k-1},\\ \underbrace{1, 1, ...,1}_{k}$, respectively.
}

\section{ Preliminary }
In this section, we will  list some concepts and Lemmas which are critical in the late proofs.

As usual,  $G = (V, E)$ is a  simple  connected graph and  $|G|$ denotes the cardinality of $V$.
For $S \subseteq V(G)$ and  $F \subseteq E(G)$,    $G[S]$ is the subgraph of $G$ induced by $S$, $G - S$ is the subgraph induced by $V(G) - S$ and $G - F$ is the subgraph of G obtained by deleting $F$.  Let $w(G - S)$ be the number of components of $G - S$ and  $S$ is a cut set if $w(G - S) \geq 2$.  For a vertex $v \in V(G)$, the neighborhood of  $v $ is the set $N(v) = N_G(v) = \{w \in V(G), vw \in E(G)\}$,  $d_G(v)$ or $d(v)$  is the degree of $v$ with $ d_{G}(v) = |N(v)|$.  
 A tree $T$ is  called a pendant tree, if $T$ has at most one  vertex shared with some cycles  in $G$. A biconnected graph is a connected graph having no cut vertices and a block is a maximal biconnected subgraph of a graph. In particular, the end block contains at most one cut vertex. Let $\lfloor {x}\rfloor $  be the largest integer that is less than or equal to $x$,  $\lceil {x} \rceil$ be the smallest integer that is greater than or equal to $x$.

By the definition of Multiplicative Zagrab index, one can easily obtain the following lemmas.

\vskip 2mm {\bf Lemma 1 } \emph{ 
For $G     \in$ $\mathcal{C}$$_n^k$ with $k \leq 1$ and $n \geq 3$, if $\prod_{1,c}(G)$ or $ \prod_{2}(G)$ attains the minimal value, then $G$ is an unicyclic graph.
}

{\bf Proof.}
For $k=0$ or $1$, by the choice of $G$, one can obtain that $G$ contains at least one cycle. Otherwise, $G$ is a tree which has at least two pendant vertices.  Assume that there exists at least two cycles in $G$,  and choose  two cycles $C_1 = x_1x_2...x_1, C_2 = y_1y_2...y_1$ and a path $P=z_1z_2...z_p$ such that $V(P) \cap V(C_1) = \{z_1\}, V(P) \cap V(C_2) = 
\{z_p\}$ and $P$ has no common vertices with  any  other cycles except $C_1, C_2$. Let $N(z_1) \cap V(C_1) = \{x_{11},x_{12}\}$ and $N(z_p) \cap V(C_1) = \{x_{p1},x_{p2}\}$, and set $G' = (G -\{x_{11}z_1,x_{p1}z_p\}) \cup \{x_{11}x_{p1}\}$, then $d_{G'}(z_1) = d(z_1) -1, d_{G'}(z_p) = d(z_p) -1$. By the definitions of $\prod_{1,c}(G)$ and $\prod_2(G)$, we have $\prod_{1,c}(G') < \prod_{1,c}(G)$ and $\prod_{2}(G') < \prod_{2}(G)$, a contradiction to the choice of $G$. Thus, Lemma 1 is ture. 
$\hfill\Box$

\vskip 2mm {\bf Lemma  2 } \emph{
Let  $G'$ be a proper subgraph of a connected graph $G$, then $\prod_{1,c}(G') < \prod_{1,c}(G), \prod_2(G') < \prod_2(G)$. In particular, for $G \in$ $\mathcal{C}$$_n^k$ with $ k \geq 2$, if $\prod_{1,c}(G)$ or $\prod_2(G)$ attains the minimal value, then $G$ is a tree. 
}

{\bf Proof.} 
Since $G'$ is a proper subgraph of $G$, by the definitions of $\prod_{1,c}(G)$ and $\prod_{2}(G)$, one can easily obtain that $\prod_{1,c}(G') < \prod_{1,c}(G)$ and $\prod_2(G') < \prod_2(G)$.   
For $k \geq 2$, we proceed to prove it by the contradiction. For $k \geq 2$, assume that $G$ is not a tree, let $C$ be a cycle of $G$ and $P_1 = u_1u_2...u_p$ and $P_2 = v_1v_2...v_q$ be two pendant paths such that $V(P_1) \cap V(C) = \{u_1\}$ and $d(v_q) =1$.  Let $w_1 \in N(u_1) \cap V(C_1) $ and $G'' = (G-\{u_1w_1, v_1v_2\}) \cup\{v_2w_1\}$, then $d_{G''}(u_1) = d(u_1) - 1, d_{G''}(v_1) = d(v_1)-1$ and $G'' \in $$\mathcal{C}$$_n^k$.
By the definitions of $\prod_{1,c}(G)$ and $\prod_2(G)$, we have   $\prod_{1,c}(G'') < \prod_{1,c}(G)$, $\prod_2(G'') < \prod_2(G)$ and Lemma 2 is true.   
$\hfill\Box$

\vskip 2mm {\bf Lemma 3 } \emph{ 
If $ \prod_{2}(G)$ attains the minimal value with $k \geq 2$, then any non-pendant vertices $u, v$ of a connected graph $G$ have the property: $|d(u) - d(v)| \leq 1$.
}

{\bf Proof.}
Since $k \geq 2$,  by Lemma 2, we have $G$ must be a tree.
On the contrary, if there are two non-pendant vertices $u, v \in V(G)$ such that 
$d(u) - d(v) \geq 2$, let $x \in N(u) -N(v)$ and $G' = (G - \{ux\} ) \cup \{vx\}$, by Fact 2 and $d_{G'}(u) = d(u) - 1, d_{G'}(v) = d(v) +1, d(v) \leq d(u) - 2 < d(u) -1$,  we have
$$\frac{\prod_{2}(G)}{\prod_2(G')} = \frac{d(v)^{d(v)}d(u)^{d(u)}}{[d(v) + 1]^{d(v)+1} [d(u) - 1]^{d(u)-1}} = \frac{[\frac{d(v)^{d(v)}}{[d(v) + 1]^{d(v)+1}}]}{[\frac{[d(u) -1]^{d(u)-1}}{d(u)^{d(u)}}]} > 1, 
$$
that is, $\prod_{2}(G') < \prod_{2}(G)$,
a contradiction with the choice of $G$. Thus, Lemma 3 is true.$\hfill\Box$

 \vskip 2mm {\bf Lemma 4 } \emph{ 
If $\prod_{1,c}(G)$ or $\prod_2(G)$ attains the maximal value, then all cycles of $G$ have length $3$ except for  at most one of them with length $4$. 
}

{\bf Proof.}  
On the contrary, 
let $C_m$ be a cycle of $G$ with $C_m = v_1v_2....v_mv_1$ and $m \geq 5$, 
$G' = (G -\{v_3v_4\}) \cup \{v_1v_3, v_1v_4\}$. 
Since $G'$ has $k$ pendant vertices, then $G' \in $ $\mathcal{C}$$_n^k$.
By the definitions of $\prod_{1,c}(G)$, $\prod_2(G)$ and 
$d_{G'}(v_1) = d(v_1)+2$, we have 
 $$\frac{\prod_{1,c}(G)}{\prod_{1,c}(G')} = \frac{d(v_1)^c}{[d(v_1) + 2]^c} < 1, \frac{\prod_2(G)}{\prod_2(G')} = \frac{d(v_1)^{d(v_1)}}{[d(v_1) + 2]^{d(v_1)+2}} < 1,$$ that is, $\prod_{1,c}(G) < \prod_{1,c}(G')$ and $\prod_2(G) < \prod_2(G')$, a contradiction with the choice of $G$. We can proceed this process until all of the cycles have length $3$ or $4$.

If there exist two cycles of length  4, say $C_1 = x_1x_2x_3x_4x_1,  C_2= y_1y_2y_3y_4y_1$ in $G$.
 Since $G$ is a cactus, then there exists a vertex $x_t \in V(C_1)$ (say $x_4$) such that there are no paths connecting 
 $x_4$ and $y_1$, $x_4$ and $y_2$ in $G - \{x_1x_4,x_3x_4\}$. Otherwise, if every vertex of $V(C)$ is either connected with 
 $y_1$ or with $y_2$  in $G - \{x_1x_4,x_3x_4\}$, then there exist a cycle that shares at  least one common edge with $C_1$, a contradiction with the definition of  cactus graph.
 Let $G^* = (G -\{x_1x_4, x_3x_4, y_1y_4\} ) \cup 
\{x_1x_3,x_4y_1,x_4y_2,y_2y_4\}$.
Since $G^*$ has $k$ pendant vertices, then $G^* \in $ $\mathcal{C}$$_n^k$.
By the definitions of $\prod_{1,c}(G)$, $\prod_2(G)$ and 
 $d_{G^*}(y_2) = d(y_2)+2$, we have
 $$\frac{\prod_{1,c}(G)}{\prod_{1,c}(G^*)} = \frac{d(y_2)^c}{[d(y_2) + 2]^c} < 1, \frac{\prod_2(G)}{\prod_2(G^*)} = \frac{d(y_2)^{d(y_2)}}{[d(y_2) + 2]^{d(y_2)+2}} < 1,$$ that is, $\prod_{1,c}(G) < \prod_{1,c}(G^*)$ and  $\prod_2(G) < \prod_2(G^*)$, a contradiction with the choice of $G$ and Lemma 4 is true.   $\hfill\Box$

 \vskip 2mm {\bf Lemma 5 } \emph{ If $\prod_{1,c}(G)$ or $\prod_2(G)$ attains the maximal value, then every dense path has length $1$ except for  at most one of them with length $2$.
}

{\bf Proof.}  
On the contrary, let $C$ be a cycle and $P = v_1v_2...v_{p-1}v_{p_i}$ with $p \geq 2$ and $j \geq 1$ be a dense path  such that $V(C) \cap V(P) = \{v_1\}$
and $d(v_{p_i}) = 1$. If $p \geq 4$, 
let $G' = G  \cup \{v_1v_{p-1}\}$. Then $G' \in G[n,k]$ and
$d_{G'}(v_1) = d(v_1) +1, d_{G'}(v_{p-1}) = d(v_{p-1}) +1$. Thus, by the definition, we have 
$\prod_{1,c}(G') > \prod_1(G)$ and $\prod_{2}(G') > \prod_2(G)$, a contradiction with the choice of $G$. We can proceed this process until $p \leq 3$, that is, all of the dense paths have the length as 1 or 2.

If there exist two such paths of length 2, say $P_1 = x_1x_2x_{3j}$, $P_2 = y_1y_2y_{3j'}$ with $x_1 \in V(C_2), y_1 \in V(C_3)$ such that $d(x_{3j}) = d(y_{3j'}) = 1$ and $j, j' \geq 1$, then let $G^* =  (G -\{y_1y_2, y_2y_{31}\}) \cup \{y_1y_{31}, x_1y_2, x_2y_2\}$.
Since $G^*$ has $k$ pendant vertices, then $G^* \in G[n,k]$.
By the definitions of $\prod_{1,c}(G)$,  $\prod_2(G)$ and  $d_{G^*}(x_1) = d(x_1) +1, d_{G^*}(x_2) = d(x_2) +1$, we have 
$\prod_{1,c}(G^*) > \prod_1(G)$ and $\prod_{2}(G^*) > \prod_2(G)$, a contradiction with the choice of $G$ and Lemma 5 is true.  $\hfill\Box$

\vskip 2mm {\bf Lemma 6 } \emph{ If $\prod_{1,c}(G)$ or $\prod_2(G)$ attains the maximal value, then 
$G$ can not have both a dense path of length $2$ and a cycle of length $4$. 
}

{\bf Proof.}  
On the contrary, let $C_1 $ be a cycle, $P = y_1y_2y_{3i}$ be a dense path  such that $V(C_1) \cap V(P) = \{y_1\}$ and  $d(y_{3i}) = 1$ for $i \geq 1$,
$C_2$ be a cycle of length $4$, say $C_2 = x_1x_2x_3x_4x_1$.
By the definition of the cactus, there exists $x_t \in V(C_2)$ (say $x_2$) such that there is no paths connecting $x_2$ and $y_1$, $x_2$ and $y_1$ in $G- \{x_1x_2,x_2x_3\}$.
 Let $G' = (G - \{x_1x_2, x_2x_3\}) \cup \{x_2y_1, x_2y_2, x_1x_3\}$, 
 then $G'$ has $k$ pendent vetices, $d_{G'}(y_1) = d(y_1) +1$ and $ d_{G'}(y_2) = d(y_2) +1$. By the definitions of $\prod_{1,c}(G)$, $\prod_2(G)$,  we have 
$\prod_{1,c}(G') > \prod_1(G)$ and $\prod_{2}(G') > \prod_2(G)$, a contradiction with the choice of $G$ and Lemma 6 is true.
 $\hfill\Box$

\vskip 2mm {\bf Lemma 7 } \emph{
Let C be a cycle of $G $ in $\mathcal{C}$$_n^k$  and $u, v \in V(C)$,  if $min\{d(u), d(v)\} > 2$, then there exist a graph $G'$ such that $\prod_2(G') > \prod_2(G)$. 
}

{\bf Proof.} 
Since $min \{d(u),d(v)\} \geq 3$, without loss of generality, let $d(u ) \geq d(v) \geq 3$,
then there exist $x \in N(v) - V(C) - N(u)$, otherwise, there will be two cycles containing at least two common vertices. Let $G' = (G -\{vx\}) \cup \{ux\}$, 
we have   $ d(u) \geq d(v)  > d(v) - 1$. By Fact 2, we have 
$$\frac{\prod_2(G)}{\prod_2(G')} = \frac{d(u)^{d(u)}d(v)^{d(v)}}{[d(u)+1]^{d(u)+1}[d(v) - 1]^{d(v) - 1}} = \frac{[\frac{d(u)^{d(u)}}{[d(u)+1]^{d(u)+1}}   ]}{[\frac{[d(v) -1]^{d(v) - 1}}{d(v)^{d(v)}}]    } < 1.$$
Thus, $\prod_2(G') > \prod_2(G)$ and Lemma 7 is true. $\hfill\Box$

 \vskip 2mm {\bf Lemma 8 } \emph{ If $\prod_{2}(G)$  attains the maximal value, then any three cycles have a common vertex.
}

{\bf Proof.} 
By the definition of the cactus, any two cycles have at most one common vertex. Now assume that there exist two disjoint cycles $C_1, C_2$ contained in $G$ such that the path $P$ connecting $C_1$ and $C_2$ is as short as possible. For convenience, let $P =u_1u_2...u_p$, $V(P) \cap V(C_1) = \{u_1\}$ and $V(P) \cap V(C_2) = \{u_p\}$.

If the path $P$ has no common edges with any other cycle(s) contained in $G$ and 
 $|E(P)| \geq 2$,  let the new graph $G'  = G \cup \{u_1u_p\}$, then  $G'  \in G[n,k]$, $d_{G'}(u_1) = d(u_1)+1$, $ d_{G'}(u_2) = d(u_2) +1$. By the definition of $\prod_2(G)$, we have   $\prod_2(G') > \prod_2(G)$. 
If  $|E(P)| = 1$, without loss of generality, let $d(u_1) \geq d(u_2)$ and $C_2 = u_2v_2v_3...u_2$, we have $v_2 \notin N(u_1)$. Otherwise, there are two cycles who have the common edge contradicted with the definition of the cactus. Let $G^* = (G - \{u_2v_2\}) \cup \{u_1v_2\}$, we have  $G^* \in G[n,k]$, 
$d_{G^*}(u_1) = d(u_1)+1$ and $d_{G^*}(u_2) = d(u) - 1$. Since $d(u_1) \geq d(u_2) > d(u_2) -1$, then 
$$
\frac{\prod_{2}(G)}{\prod_{2}(G^*)} = \frac{d(u_1)^{d(u_1)}d(u_2)^{d(u_2)}}{[d(u_1) + 1]^{d(u_1)+1}[d(u_2)-1]^{d(u_2) - 1}} = \frac{ [\frac{{d(u_1)^{d(u_1)}}} {[d(u_1) + 1]^{d(u_1)+1}}]}{[{\frac{[d(u_2)-1]^{d(u_2) - 1}}{d(u_2)^{d(u_2)}}}]} < 1,
$$
 that is, $\prod_2(G^*) > \prod_2(G)$, a contradiction with the choice of $G$.

If $P$ has some common edges with some other cycle, say $C_3$, 
 by the choice of $C_1, C_2$ and the definitions of cactus graph, we have  $\{u_1\} = C_3 \cap C_1$ and $\{u_p\} = C_3 \cap C_2$. Since $min\{d(u_1), d(u_p)\} \geq 3$, by Lemma 7, we can get that there exist $G^{**}$ such that $\prod_2(G^{**}) > \prod_2(G)$, a contradiction with the choice of $G$. 

Thus,  any two cycles of  $G$ have one common vertex. By the definition of  cactus graph, we have that
any three cycles have exactly one common vertex and Lemma 8 is true.   $\hfill\Box$

\vskip 2mm {\bf Lemma 9} \emph{ 
Let $T$ be a tree attached to a vertex $v_0$ of a cycle of  $G$, if  $ \prod_{2}(G)$ attains the maximal value, then $d(v) \leq 2$ for any $v \in V(T) - \{v_0\}$.
}

{\bf Proof.}
 Choose a graph $G$ such that  $\prod_{2}(G)$ achieves the maximal value.
On the contrary, assume that $u \in  V(T) - \{v_0\}$ is of degree $r \geq 3$ and closest to a pendant vertex. For $d(u,v_0) \geq 2$, let $G' = G \cup \{uv_0\}$, we have $G' \in G[n,k]$, $d_{G'}(u) = d(u) +1$ and $d_{G'}(v_0) = d(v_0) +1$. By the definition of  $\prod_2(G)$, we can obtain that $\prod_2(G') > \prod_2(G)$, a contradiction with the choice of $G$.
 For $d(u,v_0) = 1$, let $\{y_1,y_2,...,y_{r-2}\}$ be the r-2 neighbors of $u$ such that $d(y_i, v_0) > d(u,v_0)$,   $y$ be a neighbor of $v_0$ which belongs to a cycle $C_0$.

Since $v_0uy_1$ is a pendant path of length 2, by Lemma 6, we have that every cycle has  length  3. Let $C_0 =   v_0w_1yv_0$,
 $G'' = (G -\{uy_1\}) \cup \{v_0y_1\}$ and $G''' = (G - \{v_0y\} ) \cup \{uy\}$, then $G'', G''' \in G[n,k]$, $d_{G''}(u) = d(u) -1,  d_{G''}(v_0) = d(v_0) +1$ and $d_{G'''}(u) = d(u)+1, d_{G'''}(v_0) = d(v_0) - 1$. By the definition of  $\prod_2(G)$ and Fact 2, we can obtain
$$
\frac{\prod_{2}(G)}{\prod_{2}(G'')} = \frac{d(u)^{d(u)}d(v_0)^{d(v_0)}}{[d(u)-1]^{d(u-1)}[d(v_0)+1]^{d(v_0)+1}} = \frac{[\frac{d(v_0)^{d(v_0)}}{[d(v_0)+1]^{d(v_0)+1}}]}{[\frac{[d(u)-1]^{d(u)-1}}{d(u)^{d(u)}}]} < 1, \mbox{if}\; d(v_0) \geq d(u),
$$
$$
\frac{\prod_{2}(G)}{\prod_{2}(G''')} = \frac{d(u)^{d(u)}d(v_0)^{d(v_0)}}{[d(u)+1]^{d(u+1)}[d(v_0)-1]^{d(v_0)-1}} = \frac{[\frac{d(u)^{d(u)}}{[d(u)+1]^{d(u)+1}}]}{[\frac{[d(v_0)-1]^{d(v_0)-1}}{d(v_0)^{d(v_0)}}]} < 1, \mbox{if}\; d(v_0) < d(u),
$$
that is, $\prod_2(G'') > \prod_2(G)$ and $\prod_2(G''') > \prod_2(G)$, a contradiction with the choice of $G$.
Thus, Lemma 9 is true.
 $\hfill\Box$

 \vskip 2mm {\bf Lemma 10 } \emph{ 
If $ \prod_{2}(G)$ attains the maximal value, then all attached trees are attached to a common vertex $v_0$.
}

{\bf Proof.}
On the contrary,
suppose that there exist two trees $T_1,T_2$   attached to different vertices $v_1,v_2$  of some  cycles, say $C_1,C_2$, such that $V(C_1) \cap V(T_1) = \{v_1\}, V(C_2) \cap V(T_2) = \{v_2\}$.
 By Lemma 8, all the cycles have a common vertex $v_0$. Without loss of generality, let $v_1 \neq v_0$,  we have $d(v_0) \geq 3, d(v_1) \geq 3$.
 By Lemma 7,  there exists $G'$ such that $\prod_2(G') > \prod_2(G)$, a contradiction to the choice of $G$. Thus, Lemma 10 is true.
 $\hfill\Box$

\section{Main proofs}
In this section, we will prove the main results.
For any graph $G $ in $\mathcal{C}$$_n^k$,  if $n = 1$ or $2$, then $\prod_{1,c}(G) = \prod_2(G) = 0$ or $1$, that is, all  upper and lower bounds of  Multiplicative Zagreb indices have the same values,  respectively. Thus, all of the Theorems are true. Now we may assume that $n \geq 3$.

{\bf Proof of Theorem 1.} Choose a graph  $G$ in $\mathcal{C}$$_n^k$ such that $\prod_{1,c}(G)$ achieves the minimal value. For $k \leq 1$,  by Lemma 1, $G$ is an unicyclic graph. If $k=0$, then $G$ is a cycle, that is, the degree sequence of $G$ is $\underbrace{2, 2, ..., 2}_{n}$; If $k =1$, then $G$ has only one pendant path, that is, the degree sequence of $G$ is $3, \underbrace{2, 2, ..., 2}_{n-2}, 1$. Thus,   Theorem 1 is true.

For $k \geq 2$, by the choice of $G$ and Lemma 2, we obtain that  $G$ is a tree.  If $k =2$, then $G$ is a path, that is, the degree sequence of $G$ is $ \underbrace{2, 2, ..., 2}_{n-2}, 1,1$ and Theorem 1 is true; For $k \geq 3$, if there is a vertex $v$ with $d(v) \geq k+1$,  since $G$ is a tree, then $G$ has more than $k$ pendant vertices, a contradiction to the choice of $G$.  Thus, $d(v) \leq k$ for any $v \in V(G)$. 
Now let $v$ be the vertex with maximal degree $\Delta$, if $\Delta = k$, then $G-v$ is a set of paths. Otherwise, there exists a vertex $u \in V(G) - \{v\}$ such that $d(u) \geq 3$  and since $G$ is a tree, then $G$ contains more than $k$ pendant vertices, a contradiction  to the choice of $G$. Thus, the degree sequence of $G$ is $k, \underbrace{2, 2, ... , 2}_{n-k-1}, \underbrace{1, 1, ..., 1}_{k}$.

If $\Delta < k$, then $G$ contains at least 2 cut vertices, say $u_1, u_2, ..., u_t$, such that $G - u_i$ has at least $3$ components with $i \in [1, t]$ and $t \geq 2$. Otherwise, since $G$ is a  tree,  $G$ only contains $\Delta$ pendant vertices.
 Let $P = w_1w_2...w_s$ be a path of $G- \{u_1, u_2, ..., u_t\}$ such that $w_s \in \{u_1,u_2,...,u_t\} - \{v\}$ and $P$  contains only a unique pendant vertex $w_1$. 
Set $G' = (G -\{w_{s-1}w_s\}) \cup \{w_{s-1}v\}$, 
we have $G' \in $ $\mathcal{C}$$_n^k$, $d_{G'}(v) = d(v)+1$ and $d_{G'}(w_s) = d(w_s) - 1$. Thus, by $\Delta \geq d(w_s)> d(w_s) -1$, we have 
$$\frac{\prod_{1,c}(G)}{\prod_{1,c}(G')} = \frac{\Delta^c d(w_s)^c}{(\Delta +1)^c(d(w_s) - 1)^c} 
= \frac{\frac{\Delta^c}{(\Delta +1)^c}}{\frac{[d(w_s)-1]^c}{d(w_s)^c}} > 1,$$
that is, $\prod_{1,c} (G)$ is not  minimal, a contradiction with the choice of $G$. If the maximal degree of $G'$ is still less than $k$, then we can continue this process until $\Delta = k$, thus we can find the desired graph with the degree sequence of  $k, \underbrace{2, 2, ..., 2}_{n-k-1}, \underbrace{1, 1, ..., 1}_{k}.$
Therefore, Theorem 1 is true. $\hfill\Box$

{\bf Proof of  Theorem 2.}   Choose a graph  $G $ in $\mathcal{C}$$_n^k$ such that $\prod_{1,c}(G)$ achieves the maximal value. Let $S = \{v \in V(G), d(v) = 1\}$ and $G' = G - S$.
If $|G'|=1$,  then for $k =0$,  the degree sequence of $G$ is $0$ and for $k \neq 0$, $G$ is a star, that is, its degree sequence  is $k, \underbrace{1,1,..., 1}_{k}$.
If $|G'|=2$ and for $k = 0$, there is no such simple connected  graph;
For $k \neq 0$, by {\it "Arithmetic-Mean and Geometric-Mean inequality: $x_1x_2...x_n \leq (\frac{x_1+x_2+...+x_n}{n})^n$, the  equality holds if and only if $x_1 = x_2 =... = x_n$"}, one can obtain that
 the degree sequence of $G$ is $\lceil{\frac{k}{2}}\rceil+1, \lfloor{\frac{k}{2}}\rfloor+1,  \underbrace{1,...,1}_{k}$. 
If $|G'| =3$ and $k = 0$, by Lemma 2, we can obtain that $G$ is a cycle of length 3, that is, its degree sequence   is $2,2,2$. For $k \neq 0$, it is similar to the above proof, that is, the degree sequence of $G$ is $\lceil{\frac{k}{3}}\rceil+2,\lfloor{\frac{k}{3}}\rfloor+2, k-\lceil{\frac{k}{3}}\rceil-\lfloor{\frac{k}{3}}\rfloor+2,\underbrace{1,...,1}_{k}.$ Therefore, Theorem 2 is true.
$\hfill\Box$

{\bf  Proof of  Theorem 3.}  Choose a graph $G$ in $\mathcal{C}$$_n^k$  such that $\prod_{1,c}(G)$  achieves the maximal value. By Lemma 2 and $n - k  \geq 4$, $G$ contains some cycles. For $n - k = 4$, $G - S$ contains only one cycle $C_0$, where $S = \{v \in V(G), d(v) =1\}$. If $k = 0$, by the choice of $G$, one can obtain that   $G$ is a cycle, that is, its degree sequence  is $2,2,2,2$. If $k \neq 0$ and  $|C_0| = 4$, by adding any deleted vertex back to $G-S$, one can get a new graph $G_{01}$ with degree sequence $3, 2, 2,2,1$; If $k \neq 0$ and $|C_0| = 3$, by adding back any deleted vertex to $G-S$ such that it is adjacent to the pendant vertex in $G-S$, one can obtain a new graph $G'_{01}$. Since $G_{01}$ and $G'_{01}$ have the same degree sequences, 
 by {\it Arithmetic-Mean and Geometric-Mean inequality}, we can continue to add  any deleted vertex back to $G_{01}$  or $G'_{01}$ such that it is adjacent to the nonpendant vertex of smallest degree in $G_{01}$ or $G'_{01}$. After adding back all of the deleted vertices, we can obtain the  graph  of maximal $\prod_{1,c}$-value and Theorem 3 is true.  Thus we will consider the case when $n-k \geq 5$ below. By the choice of $G$ and Lemma 4,  $G$ contains at least two cycles.

{\small\bf Claim 1.} {\it The longest path connecting only two cycles has  length at most 1. }

{\it Proof.} 
On the contrary,  let $C_l, C_{l'}$ be two cycles and $P_1= x_1x_2...x_p$ be a path such that 
$V(C_l) \cap V(P_1) = \{x_1\}, V(C_{l'}) \cap V(P_1) = \{x_p\}$. If $p\geq 3$,  set $G' = G \cup \{x_1x_p\}$, then $d_{G'}(x_1) = d(x_1)+1$ and $d_{G'}(x_2) = d(x_2) +1$. By the definition of $\prod_{1,c}(G)$, we have
 $\prod_{1,c}(G') > \prod_{1,c}(G)$, a contradiction to the choice of $G$. Thus, $p \leq 2$ and Claim 1 is true.
$\hfill\Box$

We first deal with the case when $k = 0$.

{\small\bf Claim 2.} {\it Any three cycles have no common vertex if $k = 0$.
}

{\it Proof.} 
On the contrary, let $C_1, C_2, C_3$ be the cycles of $G$ such that $\cap _{i= 1}^3V(C_i) = \{v_0\}$, and $N(v_0) \cap V(C_i) = \{v_{i1}, v_{i2}\}$ for $i \in [1, 3]$.
 Choose $v$ of degree 2 such that $v$ is in an end block $C_t$  of $G$  and $N(v) \cap V(C_t) = \{v_{t1}, v_{t2}\}$. 
Set $G'' = (G-\{v_{21}v_0, v_{22}v_0\}) \cup \{v_{21}v, v_{22}v\}$, then $d_{G''} (v_0) = d(v_0) - 2$ and $d_{G''}(v) = d(v) + 2$. Since $d(v_0) -2  \geq 4 > d(v)$,  By the definitions of $\prod_{1,c}(G)$ and Fact 1,  we have 
$$
\frac{\prod_{1,c}(G)}{\prod_{1,c}(G'')} = \frac{d(v_0)^cd(v)^c}{[d(v_0) - 2]^c[d(v)+2]^c} = \frac{[\frac{d(v)^c}{[d(v)+2]^c}]}{[\frac{[d(v_0)-2]^c}{d(v_0)^c}]} < 1,
$$
that is, $\prod_{1,c}(G'') > \prod_1(G)$, a contradiction to the choice of $G$.
$\hfill\Box$

{\small\bf Claim 3.} {\it Every vertex  of $G$ has the degree 2, 3 or 4 if $k = 0$.
}

{\it Proof.} 
We will prove it by the contradiction. If there is a vertex $w_1$ with $d(w_1) \geq 5$, by Claim 2, we can assume that
 there are two cycles $C_4, C_{4'}$ and a path $P_2$ such that $V(C_4) \cap V(C_{4'}) \cap V(P_2) = \{w_1\}$, since $k =0$ and $G$ is a cactus, there exists a vertex $w_0$ of an end block such that $d(w_0) = 2$, that is, $d(w_0) < d(w_1) -2$. Without loss of generality,  assume that $w_0$ is closer to $C_{4'}$, let $N(w_1) \cap V(C_4) = \{w_2,w_3\}$ and $G''' = (G - \{w_1w_2, w_1w_3\}) \cup \{w_0w_2, w_0w_3\}$, by the definition of $\prod_{1,c}(G) $ and Fact 1, we have 
$$\frac{\prod_{1,c}(G)}{\prod_{1,c}(G''')} = \frac{d(w_1)^{c}d(w_0)^{c}}{[d(w_1) -2]^{c} [d(w_0) +2]^{c}} = \frac{[\frac{d(w_0)^{c}}{[d(w_0) + 2]^{c}}]}{[\frac{[d(w_1) -2]^{c}}{d(w_1)^{c}}]} < 1,
$$
that is, $\prod_{1,c}(G''') > \prod_{1,c}(G)$, a contradiction to the choice of $G$.
Thus,  Claim 3 is true. $\hfill\Box$

{\small\bf Claim 4.} {\it  There do not exist two paths of length 1 such that every path connects with only  two cycles if $k = 0$.
}

{\it Proof.}
 On the contrary, assume that there are two such paths $P_5 = z_1z_2$, $P_6 = y_1y_2$ with $z_1 \in C_6, z_2 \in C_7, y_1 \in C_8, y_2 \in C_9$ such that $N(y_1) \cap V(C_8) = \{y_{11}, y_{12}\}$ and 
$d(z_1) = d(z_2) = d(y_1) =d(y_2) = 3 $. Let $G^* =  (G- \{y_1y_2, y_1y_{11}, y_1y_{12}\}) \cup \{y_2y_{11}, y_2y_{12}, z_1y_1,z_2y_1\}$.
Since $d_{G^*}(z_1) = d_{G^*}(z_2) = d_{G^*}(y_2) = 4 , d_{G^*}(y_1) = 2$. 
By the definition of $\prod_{1,c}(G)$ and Fact 1, we have 
$$\frac{\prod_{1,c}(G)}{\prod_{1,c}(G^*)} = \frac{d(z_1)^{c}d(z_2)^{c}d(y_1)^{c}d(y_2)^{c}}{[d(z_1) +1]^{c} 
[d(z_2) +1]^{c} [d(y_1) -1]^{c} [d(y_2) +1]^{c}
} = \frac{3^c3^c3^c3^c}{4^c4^c2^c4^c} < 1,
$$
that is,  $\prod_{1,c}(G^*) > \prod_{1,c}(G)$, a contradiction to the choice of $G$ and Claim 4 is true.$\hfill\Box$

\vskip 2mm {\bf Claim 5} \emph{ 
$G$ can not have both a cycle of length $4$ and a path of length 1 connecting only with two cycles if k = 0. 
}

{\bf Proof.}  
On the contrary, let $C_{10}, C_{11}, C_{12} $ be the cycles and $P = w_1w_2$ be a path  such that $V(C_{10}) \cap V(P) = \{w_1\}$, $V(C_{11}) \cap V(P) = \{w_2\}$. If $|C_{10}| = |C_{11}| = 3$  and $|C_{12}| = 4$, then there exists a vertex $w_3 \in V(C_{12})$ such that  $d(w_3) = 3$ or $4$.  Let $C_{12} = w_3x_2x_3x_4w_3$ and $G^{**}= (G - \{w_1w_2, x_2x_3\}) \cup \{w_2w_3, w_2x_2, w_3x_3\}$, 
then $d_{G^{**}}(w_1) = d(w_1) - 1 = 2, d_{G^{**}}(w_2) = d(w_2)+1 = 4, d_{G^{**}}(w_3) = d(w_3) + 2 = 5$ or $6$ and
 $G^{**}$ has no pendent vetices. By the definitions of $\prod_{1,c}(G)$,  we have 
$$\frac{\prod_{1,c}(G)}{\prod_{1,c}(G^{**})} = \frac{d(w_1)^{c}d(w_2)^{c}d(w_3)^{c}}{[d(w_1) -1]^{c} 
[d(w_2) +1]^{c} [d(w_3) +2]^{c} 
} = \frac{3^c3^c3^c}{2^c4^c5^c} \mbox{or} \frac{3^c3^c4^c}{2^c4^c6^c} < 1,
$$
that is,  $\prod_{1,c}(G^{**}) > \prod_{1,c}(G)$,
 a contradiction with the choice of $G$. 

If $|C_{10}| = |w_1w_{12}w_{13}w_{14}w_1| =4$ and $|C_{11}|=|w_2w_{22}w_{23}w_2| = 3$, then $d(w_1) = d(w_2) = 3$, $d(w_{14}) = 2$ or $4$.  Let $G^{***} = (G- \{w_1w_{12}\}) \cup \{w_{12}w_{14}, w_{2}w_{14}\}$, we have $G^{***} \in$$\mathcal{C}$$_n^k$, $d_{G^{***}}(w_{1}) = d(w_{1})-1, d_{G^{***}}(w_{14}) = d(w_{14}) + 2, d_{G^{***}}(w_{2}) = d(w_{2}) + 1$. By the definitions of $\prod_{1,c}(G)$,  we have 
$$\frac{\prod_{1,c}(G)}{\prod_{1,c}(G^{***})} = \frac{d(w_1)^{c}d(w_{14})^{c}d(w_2)^{c}}{[d(w_1) -1]^{c} 
[d(w_{14}) +2]^{c} [d(w_2) +1]^{c} 
} = \frac{3^c2^c3^c}{2^c4^c4^c} \mbox{or} \frac{3^c4^c3^c}{2^c6^c4^c}  < 1.
$$
that is,  $\prod_{1,c}(G^{***}) > \prod_{1,c}(G)$,
 a contradiction with the choice of $G$ and Claim 5 is true. 
 $\hfill\Box$

Thus, for $k = 0$ and $n = 5$,   by the choice of $G$ and Lemma 4,  there exist two cycles of length 3, that is, its degree sequence is $4,2,2,2,2$;
For $n=6$, $G$ can be $G_l$ or $G_s$ such that $G_l$ contains two cycles of length 3 or $G_s$ contains one cycle of length 3 and one cycle of length 4, that is, the degree sequences  are $3,3,2,2,2,2$  and $4, 2,2,2,2,2$.  Since $\prod_{1,c}(G_l) > \prod_{1,c}(G_s) $,  then $\prod_{1,c}(G_l)$ attains the maximal  value;
Similarly, for $n \geq 7$, if $n = 2t+5$ with $t \geq 1$, then $G^a$ contains only the cycles of length 3 and its degree sequence  is $\underbrace{4,4,..,4}_{t+1},  \underbrace{2,2,..,2}_{t+4}$; 
If $n= 2(t+3)$, then $G^b$ contains some cycles of length 3 and one path of length 1 that connects only two cycles,  that is, its degree sequence is $\underbrace{4,4,..,4}_{t},3,3,  \underbrace{2,2,..,2}_{t+4}$.

Now we consider the case when $k \neq 0$ and define 
the following algorithm, say {\it Pro} :
 Step 1. Build $G_{T_0}$ by deleting all the dense paths such that $G_{T_0}$ satisfies the case of $k = 0$, that is, $G_{T_0}$ is either $G^a$ or $G^b$;
 Step 2.  Build $G_{T_i}$ by adding a deleted path to $G_{T_{i-1}}$  such that   
it  is adjacent to a non-pendant vertex of smallest  degree in $G_{T_{i-1}}$, $i \geq 1$; 
Step 3. Stop, if there is no remaining deleted  paths; Go to Step 2, if otherwise. 

By the choice of $G$ and Lemma 5,  all of the dense paths of $G$ have length  1 except for at most one of them with length 2. If all of the dense paths of $G$ have  length 1,  by {\it Arithmetic-Mean and Geometric-Mean inequality}, we can directly use {\it Pro} to get a new graph $G_T$ of maximal $\prod_{1,c}$-value. Thus,  for $k <  4+t$, $G_T$ contains no cycles of length greater than 3, no dense paths of length greater than 1,  no paths of length greater 0 that connects only two cycles except for at most one of them with length 1 and $d_{G_T}(w_a) \in \{2,3,4\}$, where  $w_a$ is any  nonpendant vertex of $ G_T$;  For $k \geq 4+t$, we have $|d_{G_T}(w_b) - d_{G_T}(w_c)| \leq 1$,  where $w_b, w_c$ are any  nonpendant vertices of $ G_T$, that is, the statement $(i)$ or $(ii)$ is true.   
If there is one of the dense paths of $G$ with length 2,  then  set  $P_1 = u_1u_2u_{31}$, $P_2 = u_1u_2u_{32},...,P_{r-1} = u_1u_2u_{3(r-1)}$ with $ d(u_{3i}) =1, i \in [1, r-1]$. 
 By {\it Arithmetic-Mean and Geometric-Mean inequality}, we can  use {\it Pro} to get a new graph $G_{T}$ such that  $\prod_{1,c}(G_T) \geq \prod_{1,c}(G)$.  

By the proof of the case for $k = 0$,
if $G_T$ contains a  path $P_T = w_{T1}w_{T2}$ connecting only two cycles, say $C_{T1}, C_{T2}$, such that $V(C_{T1}) \cap V(P_T) = \{w_{T1}\}, V(C_{T2}) \cap V(P_T) = \{w_{T2}\}$, then set $G_1 =  (G_{T} - \{u_2u_{i}, i \in [1, r-1]  \}) \cup \{u_{31}w_{T1}, u_{31}w_{T2}, u_{31}u_j, i \in[2, r-1]\}$. Since $G_1 \in$ $\mathcal{C}$$_n^k$, $d_{G_1}(w_{T1}) = d(w_{T1})+1, d_{G_1}(w_{T2}) = d(w_{T2}) + 1, d(u_2) = r, d_{G_1}(u_2) = 1,  d(u_{31}) =1$ and $d_{G_1}(u_{31}) = r$,  by the definition of $\prod_{1,c}(G)$ and Fact 1, we have 
$$
\frac{\prod_{1,c}(G)}{\prod_{1,c}(G_1)} = \frac{d(w_{T1})^{c}d(w_{T2})^{c}d(u_{31})^{c}d(u_2)}{[d(w_{T1})+1]^{c}[d(w_{T2})+1]^{c}d(u_{31})^{c}d_1(u_2)}  < 1,
$$
that is, $\prod_{1,c}(G_1) > \prod_{1,c}(G) $, a contradiction to the choice of $G$.

If $G_T$ contains no such path $P_T$ and $|d(u_2) - d(v_T)| \leq 1$ for any nonpendant vertices $v_T, v'_T$ of $G_T$, when $|d(v_T)-d(v'_T)| \leq 1$, then the statement $(i)$  is true; When there exist  $v_T, v'_T$ such that $|d(v_T)-d(v'_T)| > 1$, by the construction of $G_T$, we have
$d(v_T),d(v'_T) \in  \{2,3,4\}$ and $G$ contains  no dense paths of length greater than 1 except for at most one of them with length 2, that is, the statement (ii) is true.  Otherwise, if  there exists a vertex $v_T$ such that $|d(u_2) - d(v_T)| > 1$, then without loss of generality, choose $C_{T3}$ and $C_{T4}$ such that $V(C_{T3}) \cap V(C_{T4}) = \{v_{T}\}$ and let $N(v_{T}) = \{v_{ci}, i \geq 4\}$ such that $v_{c1}, v_{c2} \in V(C_{T3})$, $v_{c3}, v_{c4} \in V(C_{T3})$.
When $d(v_T) - d(u_2) > 1$, set $G_2 = (G_{T} - \{v_{T}v_{c1}, v_{T}v_{c2}, u_2u_{3i}, i \geq 1\}) \cup \{u_{31}v_{c1},u_{31}v_{c2}, u_{31}v_{T},   u_{31}u_{3i}, i \geq 2 \}$, then $d_{G_2}(u_2) =1, d_{G_2}(u_{31}) = d(u_2) +1$ and $d_{G_2}(v_T) = d(v_T) - 1$.  When $d(u_2) - d(v_T) > 1$, that is, $d(u_2) > 3$, set  $G_3 = (G_{T} - \{v_{T}v_{c1}, v_{T}v_{c2}, u_2u_{3i}, i \geq 1\}) \cup \{u_{31}v_{c1},u_{31}v_{c2}, u_{31}v_{T}, u_{32}v_T, u_{33}v_T,  u_{31}u_{3i}, i \geq 4 \}$, then $d_{G_3}(u_2) =1, d_{G_3}(u_{31}) = d(u_2) -1$ and $d_{G_3}(v_T) = d(v_T) + 1$. 
 Since $G_2, G_3 \in$ $\mathcal{C}$$_n^k$,  by the definition of $\prod_{1,c}(G)$ and Fact 1, we have 
$$
 \frac{\prod_{1,c}(G)}{\prod_{1,c}(G_2)} = \frac{d(u_{31})^{c}d(u_2)^cd(v_T)^c}{[d(u_2)+1]^{c}1^c[d_2(v_T) -1]^c}  < 1, \frac{\prod_{1,c}(G)}{\prod_{1,c}(G_3)} = \frac{d(u_{31})^{c}d(u_2)^cd(v_T)^c}{[d(u_2) -1]^{c}1^c[d(v_T)+1]^c}  < 1, 
$$
that is, $\prod_{1,c}(G_2) > \prod_{1,c}(G) $, $\prod_{1,c}(G_3)  >  \prod_{1,c}(G) $, a contradition to the choice of $G$.
Therefore,   Theorem 3 is true. 
$\hfill\Box$

{\bf Proof of  Theorem 4.} 
Choose a graph $G $ in $\mathcal{C}$$_n^k$ such that $\prod_2(G)$  achieves the minimal value.
 By Lemma 1, $G$ is an unicyclic graph for $k \leq 1$. If $k = 0$, then $G$ is a cycle, that is, its degree sequence  is $\underbrace{2, 2, ..., 2}_{n}$; If $k =1$, then $G$ has only one pendant path, that is, its degree sequence  is $3, \underbrace{2, 2, ..., 2}_{n-2}, 1$.

For $k \geq 2$, by Lemma 2, we only need to consider $G$  as a tree. Since $\sum_{v \in V(G)} d(v) = 2(n-1)$, then the average degree of $G$ except the pendant vertices is  $\frac{\sum_{v \in V(G)} d(v) - k }{n-k} = \frac{2(n-1) - k }{n- k} = 2 + \frac{k-2}{n-k} = 2 + \gamma$. By Lemma 3, if all of nonpendant vertices have  degree   $2+ \lfloor{\gamma}\rfloor$ or $2+ \lceil{\gamma}\rceil$,  then $\prod_{2}(G)$ attains the minimal value.
Set the number of the vertices with degree $2+ \lfloor{\gamma}\rfloor$ to be  $y_1$,  the number of the vertices with degree $2+ \lceil{\gamma}\rceil$ to be  $y_2$, we have $y_1+y_2+k = n$ and $(2+ \lfloor{\gamma}\rfloor)y_1+(2+ \lceil{\gamma}\rceil)y_2+k = 2(n-1)$. If $\lfloor{\gamma}\rfloor = \lceil{\gamma}\rceil$, then  Theorem 4 is true; If $  \lceil{\gamma}\rceil - \lfloor{\gamma}\rfloor = 1$,
by solving the above equations, we have $y_1 = n-2k+2+ \lfloor{\gamma}\rfloor(n-k), y_2 = k-2 -\lfloor{\gamma}\rfloor (n-k), $ that is, 
 its degree sequence is
$\underbrace{2+ \lceil{\gamma}\rceil, 2+ \lceil{\gamma}\rceil, ..., 2+ \lceil{\gamma}\rceil}_{k-2 -\lfloor{\gamma}\rfloor (n-k) }, \underbrace{2+\lfloor{\gamma}\rfloor, 2+\lfloor{\gamma}\rfloor, ..., 2+ \lfloor{\gamma}\rfloor}_{n-2k+2+ \lfloor{\gamma}\rfloor(n-k)}, \underbrace{1, 1, ..., 1}_{k}$. Therefore, Theorem 4 is true.$\hfill\Box$

{\bf Proof of  Theorem 5.} 
Choose $G  $ in $\mathcal{C}$$_n^k$ such that $\prod_2(G)$ achieves the maximal value.
 By Lemma 4, the lengths of all cycles in $G$ are $3$ except for at most one of them with length  $4$; By Lemmas 5 and 9, every pendant path has length of $1$  except for at most one of them with length  $2$. By Lemma 6, $G$ can not have both a dense path of length $2$ and a cycle of length $4$; By Lemma 8, any three cycles have a common vertex $v_0$; By Lemma 10, any tree attachs to the same vertex $u$. Now we  show that $u = v_0$.  Otherwise, if $u \neq v_0$ and  $d(v_0) \geq d(u)$, let $C^*$ be the cycle that contains $u$ and  $G' = (G-\{uy | y \in N(u) - V(C^*) \}) \cup \{v_0y | N(u) - V(C^*) \}$ with $|N(u)- V(C^*)| = t_1$, by Fact 2 and $d_{G'}(u) = d(u) - t_1, d_{G'}(v_0) = d(v_0) + t_1$, we have
$$\frac{\prod_{2}(G)}{\prod_2(G')} = \frac{d(v_0)^{d(v_0)}d(u)^{d(u)}}{[d(v_0) + t_1]^{d(v_0)+t_1} [d(u) - t_1]^{d(u)-t_1}} = \frac{\frac{d(v_0)^{d(v_0)}}{[d(v_0) + t_1]^{d(v_0)+t_1}}}{\frac{[d(u) - t_1]^{d(u)-t_1}}{d(u)^{d(u)}}} < 1, 
$$
that is, $\prod_{2}(G') > \prod_{2}(G)$,
a contradiction with the choice of $G$; If $d(u) > d(v_0)$, let
$G'' = (G - \{v_0y | y \in N(v_0) - V(C^*)\} ) \cup \{uy | y \in N(v_0) - V(C^*)\}$ with $|N(u)- V(C^*)| = t_2$,  by Fact 2 and $d_{G''}(v_0) = d(v_0) - t_2, d_{G''}(u) = d(u) + t_2$, we have
$$\frac{\prod_{2}(G)}{\prod_2(G'')} = \frac{d(u)^{d(u)}d(v_0)^{d(v_0)}}{[d(u) + t_2]^{d(u)+t_2} [d(v_0) - t_2]^{d(v_0)-t_2}} = \frac{\frac{d(u)^{d(u)}}{[d(u) + t_2]^{d(u)+t_2}}}{\frac{[d(v_0) - t_2]^{d(v_0)-t_2}}{d(v_0)^{d(v_0)}}} < 1, 
$$
that is, $\prod_{2}(G'') > \prod_{2}(G)$,
a contradiction with the choice of $G$.  Therefore, we can obtain the construction of $G$ as follows:
If $n-k \equiv 0 (\mbox{mod}\; 2)$, then the degree sequence of $G$ is 
$n-2, \underbrace{2, 2, ..., 2}_{n-k-1 }, \underbrace{1, 1, ...,1}_{k}$;
if $n-k \equiv 1 (\mbox{mod}\; 2)$, then the degree sequence of $G$ is 
$n-1, \underbrace{2, 2, ..., 2}_{n-k-1}, \underbrace{1, 1, ...,1}_{k}$. Thus,   Theorem 5 is true.
$\hfill\Box$

\section{Compliance with Ethical Standards}

{\it Conflict of interest:} Shaohui Wang and Bing Wei state that there are no conflicts of interest.
Patients rights and animal protection statements: This article does not contain any studies
with human or animal subjects.

{\it Funding:} The first author was partially supported by the Summer Graduate Research Assistantship
Program of Graduate School, the second author was partially supported by College
of Liberal Arts Summer Research Grant.

\end{document}